\documentclass[10pt,twocolumn,twoside,a4paper]{IEEEtran}

\usepackage[square,comma,numbers,sort]{natbib}
\makeatletter
\def\NAT@def@citea{\def\@citea{\NAT@separator}}
\makeatother

\usepackage{amsmath}
\allowdisplaybreaks[1]
\usepackage{amssymb}

\usepackage{txfonts}
\usepackage{nccmath}
\usepackage{nicefrac}

\usepackage{mathtools}
\newcommand\scaleddot{\scalebox{.89}{.}}

\makeatletter
\renewcommand{\dddot}[1]{%
  {\mathop{\kern\z@#1}\limits^{\makebox[0pt][c]{\vbox to-2.2\ex@{\kern-\tw@\ex@
   \hbox{\normalfont\scaleddot\kern-0.5pt\scaleddot\kern-0.5pt\scaleddot}\vss}}}}}
\renewcommand{\ddddot}[1]{%
  {\mathop{\kern\z@#1}\limits^{\makebox[0pt][c]{\vbox to-2.2\ex@{\kern-\tw@\ex@
   \hbox{\normalfont\scaleddot\kern-0.5pt\scaleddot\kern-0.5pt\scaleddot\kern-0.5pt\scaleddot}\vss}}}}}
\makeatother

\newcommand{\invsign}{\scalebox{0.5}[0.75]{\( - \)}1}
\newcommand{\invtrans}{\scalebox{0.5}[0.75]{\( - \)}\top}

\renewcommand{\S}{\mathbf{S}}
\newcommand{\SA}{\mathbf{A}}
\newcommand{\SB}{\mathbf{B}}
\newcommand{\SC}{\mathbf{U}}
\newcommand{\SD}{\mathbf{V}}

\newcommand{\T}{\mathbf{T}}
\newcommand{\TA}{\mathbf{C}}
\newcommand{\TB}{\mathbf{D}}
\newcommand{\TC}{\mathbf{W}}
\newcommand{\TD}{\mathbf{Z}}

\newcommand{\x}{\mathbf{x}}

\newcommand{\y}{\mathbf{y}}

\newcommand{\w}{\mathbf{w}}
\renewcommand{\v}{\mathbf{v}}
\newcommand{\s}{\mathbf{s}}
\newcommand{\Y}{Y}
\newcommand{\m}{\mathbf{m}}

\renewcommand{\r}{\mathbf{r}}

\renewcommand{\d}{\mathrm{d}}

\renewcommand{\a}{\mathbf{a}}
\newcommand{\f}{\mathbf{f}}
\newcommand{\h}{\mathbf{h}}

\renewcommand{\P}{\mathbf{P}}

\newcommand{\M}{\mathbf{M}}

\newcommand{\N}{\mathbf{N}}

\newcommand{\F}{\mathbf{F}}
\newcommand{\G}{\mathbf{G}}
\renewcommand{\H}{\mathbf{H}}
\newcommand{\I}{\mathbf{I}}
\newcommand{\J}{\mathbf{J}}
\renewcommand{\O}{\mathbf{0}}

\begin{document}
\title{Marginalized Bayesian filtering with Gaussian priors and posteriors}

\author{John-Olof~Nilsson
\thanks{John-Olof~Nilsson is with the Dept. of Signal Processing, KTH Royal Institute of Technology, Osquldas v. 10, $100\,44$ Stockholm, Sweden. jnil02@kth.se\protect}}

\maketitle

\begin{abstract}
Marginalization techniques are presented for the Bayesian filtering problem under the assumption of Gaussian priors and posteriors and a set of sequentially more constraining state space model assumptions. The techniques provide the capabilities to 1) reduce the filtering problem to essential marginal moment integrals, 2) combine model and numerical approximations to the moment integrals, and 3) decouple modelling and system composition. Closed-form expressions of the posterior means and covariances are developed as functions of subspace projection matrices, subsystem models, and the marginal moment integrals. Finally, we review prior work and how the results relate to Kalman and marginalized particle filtering techniques.
\end{abstract}

\section{Introduction}
The Bayesian filtering problem is in general intractable and direct numerical approximations, such as particle filters, are only applicable to low-dimensional systems. Beyond particle filters, marginalization, combined with a (conditionally) static, linear, and Gaussian state space subspace, can confine the intractable part to the complementary subspace making particle filter techniques, so-called marginalized or Rao-Blackwellized particle filters, applicable to higher dimensional systems. Unfortunately, for commonly appearing dynamic subspace structures, one is in general referred to filtering with the assumptions of Gaussian priors and posteriors, e.g. Kalman filtering. However, marginalization may still be highly valuable and used to exploit structures in the models. In particular, in this article I present general marginalization techniques which give the following model structure exploiting capabilities:


\emph{1) Focusing numerical approximations:} The Gaussian prior and posterior assumptions reduce the Bayesian filtering to a set of intractable moment integrals. Consequently, the related filtering techniques are essentially defined by the different ways they approximate the moment integrals. However, with subspace structures, only subdimensions of the integrals are intractable and marginalization can reduce the integrals to those intractable dimensions, thereby focusing the numerical approximation techniques.

\emph{2) Combining model and numerical approximations:} Traditionally, nonlinear non-Gaussian models are handled by linearization and Gaussian approximations, giving different flavors of the Kalman filters. However, linearizing or making a Gaussian approximation of a part of a model is just a way of introducing related linear and Gaussian subspaces. Then, marginalization can be used to eliminate such subspaces from the moment integrals, hence enabling combination of model and numerical approximations.

\emph{3) Decoupling modelling and system composition:} Filtering problems often include multiple information sources and hence multiple models. The different output and state models are typically dependent on different system states, creating dynamic subspace structures. However, the specific subspace structures will be dependent on the system composition, creating an implicit coupling between the modelling and the composition. In this case, marginalization can be used to decouple the subsystem modelling and the system composition such that they can be performed separately.

The techniques give extended Bayesian filtering abilities. The novelty is both within the capabilities themselves and the fact that they are derived under different (weaker) assumptions than comparative filtering techniques. In combination with moment integral approximations, the results constitute a new class of Bayesian filtering techniques.
The core results of this article are closed form expressions of the posterior means and covariances in terms of subspace projection matrices, subsystem models, and explicit marginal moment integrals, for a set of sequentially more constraining state space model assumptions. 
In the following subsection, we pose the Bayesian filtering problem and the state space model assumptions. The expressions for the posterior means and covariances are developed over Sections~\ref{sec:inactive_subspace_marginalization}-\ref{sec:Closed_form_solutions_for_the_linear_case}. Following this, in Sections~\ref{sec:Marginal_moments_approximation}-\ref{sec:system_composition} we return to how they give capabilities 1-3. Finally, in Section~\ref{sec:connections} we review prior work and applications and how the results connect to related filtering techniques.
Due to the wide range of analytical results and the more specific nature of suitable moment integral approximations, we stay away from specific approximation techniques and simulations.

\vspace{1.25mm}



\emph{Notation:} Let $\mathbf{A}^{\bot\top}=\mathbf{A}\mathbf{A}^\top$ and $\mathbf{A}^{\bot+\top}=\mathbf{A}+\mathbf{A}^\top$. $\mathbf{0}$ and $\I$ are zero and identity matrices of appropriate dimensions. $[\cdot\:\cdot]$ and $[\cdot\,;\cdot]$ are row and column vectors, respectively. $p(\a|\mathbf{b})$ is the probability density function (pdf) of $\a$ given (where applicable) $\mathbf{b}$. 
$\mathcal{N}_\a(\mathbf{c},\P)$ is the Gaussian pdf of $\a$ with mean $\mathbf{c}$ and covariance $\P$. $\delta(\cdot)$ is the Dirac delta function.
Conditioning on output up to instant $k$ is denoted with $\hat{\a}$ and $k-1$ with $\check{\a}$.

\subsection{Problem formulation}
Consider the general discrete state space system
\begin{equation*}
\begin{matrix}
\x_k=f_k(\x_{k-1},\w_k)\\
\y_k=h_k(\x_{k-1},\v_k)
\end{matrix}
\quad\Rightarrow\quad
\begin{matrix}
p(\x_k|\x_{k-1})\\
p(\y_k|\x_k)
\end{matrix}
\end{equation*}
where $\x_k\in\mathbb{R}^n$ and $\y_k\in\mathbb{R}^m$ and where $\w_k\sim p(\w_k)$ and $\v_k\sim p(\v_k)$ are independent. In the latter implied (and sufficient) density representation, $p(\x_k|\x_{k-1})$ and $p(\y_k|\x_k)$ captures both $f_k(\cdot,\cdot)$ and $h_k(\cdot,\cdot)$ and $p(\w_k)$ and $p(\v_k)$. For $p(\y_k|\x_k)$ we do not require quantified output but it may also encode any side information $\y_k$, e.g. constraints.
The problem is to recursively calculate the conditional densities $p(\x_k|\Y_{k-1})$ and $p(\x_k|\Y_k)$ given $p(\x_0|\Y_{0})$, where $\Y_k=\{\y_k\}$ is the set of all output up to instant $k$. This is the Bayesian filtering problem which in general cannot be solved analytically. However, the intractability can be limited to specific subspaces with further assumptions:

\vspace{2mm}
\noindent$\mathcal{N}$\emph{) Gaussian priors and posteriors:} While $p(\x_k|\x_{k-1})$ and $p(\y_k|\x_k)$ are arbitrary densities, the conditional state densities, 
\begin{equation*}
p(\x_k|\Y_{k-1})=\mathcal{N}_{\x_k}(\check\x_k,\check\P_k)\:\:\:\text{and}\:\:\:
p(\x_k|\Y_{k})=\mathcal{N}_{\x_k}(\hat\x_k,\hat\P_k)
\end{equation*}
i.e. the initial and estimated distributions, are Gaussian.

\vspace{2mm}
\noindent\emph{a) Dynamic active subspaces:} There are matrices $\SA_k\in\mathbb{R}^{\alpha_k,n}$, $\SB_k\in\mathbb{R}^{n-\alpha_k,n}$, $\TA_k\in\mathbb{R}^{\beta_k,n}$, and $\TB_k\in\mathbb{R}^{n-\beta_k,n}$, and corresponding functions $f^a_k(\cdot,\cdot)$ and $h^a_k(\cdot,\cdot)$, such that
\begin{align*}
\S_k\x_{k}&=\left[f^a_k(\SA_k\x_{k-1},\w_k);\SB_k\x_{k-1}\right]\\
\y_k&=h^a_k(\TA_k\x_k,\v_k)
\end{align*}
where $\S_k=[\SA_k;\SB_k]$ and $\T_k=[\TA_k;\TB_k]$ are invertible.

\vspace{2mm}
\noindent The structure of assumption \emph{a} frequently arises due to different models being functions of different state space subspaces $\SA_k\x_{k-1}$ and $\TA_k\x_k$.
$\mathcal{N}$ and \emph{a} are the basic assumptions here. Results for the filtering problem will be presented for this general case. However, further structure is often available, structure which can greatly simplify the filtering. Consequently, results will also be presented for three further assumptions about the active subspaces:

\vspace{2mm}
\noindent\emph{b) Conditionally additive noise:}
There are functions $f^b_k(\cdot)$ and $h^b_k(\cdot)$, and matrices $\G^b_k(\SA_k\x_k)$ and $\J^b_k(\TA_k\x_k)$ such that
\begin{align*}
\SA_k\x_{k}&=f^b_k(\SA_k\x_{k-1})+\G^b_k(\SA_k\x_{k-1})\w_k\\
       \y_k&=h^b_k(\TA_k\x_k)+\J^b_k(\TA\x_k)\v_k
\end{align*}
where $p(\w_k)$ has zero mean and covariance $\P_{\w_k}$ and $\J^b_k(\TA\x_k)$ is invertible. (The mean and covariance of $\v_k$ may be undefined.)

\vspace{2mm}
\noindent\emph{c) Conditionally linear and Gaussian subspaces:}
There are matrices $\SA_k=[\SA^n_k;\SA^l_k]$, $\TA_k=[\TA^n_k;\TA^l_k]$, $\F^c_k(\SA^n_k\x_k)$, $\G^c_k(\SA^n_k\x_k)$, $\H^c_k(\TA^n_k\x_k)$, and $\J^c_k(\TA^n_k\x_k)$ and functions $f^c_k(\cdot)$ and $h^c_k(\cdot)$ such that
\begin{align*}
\SA_k\x_{k}&=f^c_k(\SA^n_k\x_{k-1})+\F^c_k(\SA^n_k\x_{k-1})\SA^l_k\x_{k-1}\!+\!\G^c_k(\SA^n_k\x_{k-1})\w_k\\
       \y_k&=h^c_k(\TA^n_k\x_k)+\H^c_k(\TA^n\x_k)\TA^l_k\x_k+\J^c_k(\TA^n\x_k)\v_k
\end{align*}
where $\w_k\sim\mathcal{N}_{\w_k}(\mathbf{0},\P_{\w_k})$ and $\v_k\sim\mathcal{N}_{\w_v}(\mathbf{0},\P_{\v_k})$, i.e. there are conditionally linear subspaces with additive Gaussian noise within the active subspaces.

\vspace{2mm}
\noindent\emph{d) Dynamic active linear and Gaussian subspaces:}
There are matrices $\F^d_k$, $\G^d_k$, $\H^d_k$ and $\J^d_k$, and vectors $\f^d_k$ and $\h^d_k$ such that
\begin{align*}
\SA_k\x_{k}&=\f^d_k+\F^d_k\SA_k\x_{k-1}+\G^d_k\w_k\\
       \y_k&=\h^d_k+\H^d_k\TA_k\x_k+\J^d_k\v_k
\end{align*}
where $p(\w_k)$ has zero mean and covariance $\P_{\w_k}$ and $\v_k\sim\mathcal{N}_{\v_k}(\mathbf{0},\P_{\v_k})$, i.e. the subspace state equation and the output equation are linear and Gaussian.

\vspace{2mm}
\noindent Note that the assumptions about the state equation and the output equation within the above assumptions are independent. Consequently, they may be used separately. However, since the analysis for the prediction and the update are similar, they are treated jointly.
Further, $p(\w_k)$ having a defined mean and covariance $\P_{\w_k}$ is strictly not required for the results of assumption \emph{a}. However, it is in general required for assumption $\mathcal{N}$ to make sense. Similarly, $\w_k\sim\mathcal{N}_{\w_k}(\mathbf{0},\P_{\w_k})$ is not strictly required for the results for assumption \emph{d} but it is implicit from the linear state equation and assumption $\mathcal{N}$. Hence, the assumptions \emph{a-d} are sequentially more constraining.

\section{Inactive subspace marginalization}\label{sec:inactive_subspace_marginalization}
Assumption $\mathcal{N}$ implies that only the means, $\check\x_k$ and $\hat\x_k$, and the covariances, $\check\P_k$ and $\hat\P_k$, need to be calculated, instead of arbitrary densities $p(\x_k|\Y_{k-1})$ and $p(\x_k|\Y_{k})$.
If, additionally, assumption \emph{a} is made, the desired means and covariances can be expressed as functions of marginal moments in the nonlinear non-Gaussian subspaces, as will be shown in the following subsections.

\vspace{4mm}
\emph{For brevity, hereinafter, indices of $f_k(\cdot,\cdot)$, $h_k(\cdot,\cdot)$, $\S_k$ and $\T_k$ and related quantities are dropped, while retaining implicit dependencies on $k$. Further, required indices of states, covariances, and outputs are indicated with diacritic dots, $\a_k=\dot{\a}$ and $\a_{k-1}=\ddot{\a}$.}

\subsection{Marginalization for prediction}
Let $\SC\in\mathbb{R}^{n,a}$ and $\SD\in\mathbb{R}^{n,n-a}$ such that $\S^{\invsign}=[\SC\:\SD]$. Further, let $\M=\SB\check{\ddot{\P}}\SA^\top(\SA\check{\ddot{\P}}\SA^\top)^{\invsign}$. Assumption \emph{a} implies $p(\S\dot\x|\ddot\x)=p(\SA\dot\x|\SA\ddot\x)\delta(\SB\dot\x-\SB\ddot\x)$. This, and assumptions $\mathcal{N}$ and \emph{a}, imply
\begin{gather*}
\int{p}(\S\dot\x,\S\ddot\x|\ddot\Y)\d(\SB\ddot\x)\d(\SB\dot\x)=p(\SA\dot\x,\SA\ddot\x|\ddot\Y)\\ \int{p}(\S\dot\x,\S\ddot\x|\ddot\Y)\d(\SA\ddot\x)\d(\SA\dot\x)=\delta(\SB\dot\x-\SB\ddot\x)p(\SB\ddot\x|\ddot\Y)\\
\int\SB\ddot{\x}p(\S\ddot{\x}|\ddot{\Y})\d(\SB\ddot{x})=(\SB\check{\ddot{\x}}+\M(\SA\ddot{\x}-\SA\check{\ddot\x}))p(\SA\ddot{\x}|\ddot{\Y}).
\end{gather*}
Finally, introduce the marginal moments
\begin{gather*}
\SA\check{\dot\x}=\int\SA\dot{\x}p(\SA\dot{\x},\SA\ddot{\x}|\ddot{\Y})\d(\SA\ddot{\x})\d(\SA\dot{\x})\\
\check{\P}_{\SA\dot\x}=\int(\SA\dot\x-\SA\check{\dot\x})^{\bot\top}p(\SA\dot\x,\SA\ddot\x|\ddot\Y)\d(\SA\ddot{\x})\d(\SA\dot{\x})\\
\check{\P}_{\SA\dot\x,\SA\ddot\x}=\int(\SA\dot\x-\SA\check{\dot\x})(\SA\ddot\x-\SA\check{\ddot\x})^\top{p}(\SA\dot{\x},\SA\ddot{\x}|\ddot{\Y})\d(\SA\ddot{\x})\d(\SA\dot{\x}).
\end{gather*}
Together, this gives the predicted mean
\begin{align*}
\check{\dot{\x}}&=\int\dot{\x}p(\dot{\x}|\ddot\x)p(\ddot\x|\ddot\Y)\d\ddot\x\d\dot\x\\
          &=\S^{\invsign}\int\S\dot\x{p}(\S\dot\x|\ddot\x)p(\S\ddot\x|\ddot{\Y})\d(\S\ddot\x)\d(\S\dot\x)\\
          &=\S^{\invsign}\int[\SA\dot\x;\SB\dot\x]p(\SA\dot{\x}|\SA\ddot{\x})\delta(\SB\dot\x-\SB\ddot\x)p(\S\ddot{\x}|\ddot{\Y})\d(\S\dot{\x})\d(\S\ddot{\x})\\
          &=\SC\int\SA\dot{\x}p(\SA\dot{\x}|\SA\ddot{\x})p(\SA\ddot{\x}|\ddot{\Y})\d(\SA\dot{\x})\d(\SA\ddot{\x})+\SD\SB\check{\ddot{\x}}\\
          &=\SC\SA\check{\dot\x}+\SD\SB\check{\ddot{\x}}.
\end{align*}
Further, note that $\SB\check{\dot\x}=\SB\SC\SA\check{\dot\x}+\SB\SD\SB\check{\ddot\x}=\SB\check{\ddot\x}$. Together, this gives the predicted covariance
\begin{align*}
\check{\dot{\P}}&=\int(\dot\x-\check{\dot\x})^{\bot\top}p(\dot\x,\ddot\x|\ddot\Y)\d\ddot\x\d\dot\x\\
                &=\int\S^{\invsign}\S(\dot\x-\check{\dot\x})^{\bot\top}\S^\top\S^{\invtrans}p(\S\dot\x,\S\ddot\x|\ddot\Y)\d(\S\ddot\x)\d(\S\dot\x)\\
                &=\int\Big(\SC(\SA\dot\x-\SA\check{\dot\x})^{\bot\top}\SC^\top+\SD(\SB\dot\x-\SB\check{\dot\x})^{\top\bot}\SD^\top\\ &\quad\quad\quad\quad+(\SC\SA(\dot\x-\check{\dot\x})^{\bot\top}\SB^\top\SD^\top)^{\bot+\top}\Big){p}(\S\dot\x,\S\ddot\x|\ddot\Y)\d(\S\ddot\x)\d(\S\dot\x)\\
                &=\int\SC(\SA\dot\x-\SA\check{\dot\x})^{\bot\top}\SC^\top{p}(\SA\dot\x,\SA\ddot\x|\ddot\Y)\d(\SA\ddot\x)\d(\SA\dot\x)\\
                &\quad+\int\SD(\SB\dot\x-\SB\check{\dot\x})^{\bot\top}\SD^\top\delta(\SB\dot\x-\SB\ddot\x)p(\SB\ddot\x|\ddot\Y)\d(\SB\ddot\x)\d(\SB\dot\x)\\
                &\quad+\int\Big(\SC(\SA\dot\x-\SA\check{\dot\x})(\SB\dot\x-\SB\check{\dot\x})^\top\SD^\top\Big)^{\bot+\top}\\
                &\quad\quad\quad\quad\quad\quad\quad\quad\:{p}(\SA\dot{\x}|\SA\ddot{\x})\delta(\SB\dot\x-\SB\ddot\x)p(\S\ddot{\x}|\ddot{\Y})\d(\S\ddot{\x})\d(\S\dot{\x})\\
                &=\SC\check{\P}_{\SA\dot\x}\SC^\top+\SD\SB\check{\ddot\P}\SB^\top\SD^\top\\
                &\quad+\int\Big(\SC(\SA\dot\x-\SA\check{\dot\x})(\SA\ddot\x-\SA\check{\ddot\x})^\top\M^\top\SD^\top\Big)^{\bot+\top}\\
                &\quad\quad\quad\quad\quad\quad\quad\quad\quad\quad\quad\quad\:{p}(\SA\dot{\x}|\SA\ddot{\x})p(\SA\ddot\x|\ddot\Y)\d(\SA\ddot{\x})\d(\SA\dot{\x})\\
                &=\SC\check{\P}_{\SA\dot\x}\SC^\top+\SD\SB\check{\ddot\P}\SB^\top\SD^\top+(\SC\check{\P}_{\SA\dot\x,\SA\ddot\x}\M^\top\SD^\top)^{\bot+\top}.
\end{align*}
Consequently, with assumptions $\mathcal{N}$ and \emph{a}, obtaining $p(\dot{\x}|\ddot{\Y})$ from $p(\ddot{\x}|\ddot{\Y})$, i.e. the prediction, boils down to calculating the marginal moments $\SA\check{\dot\x}$, $\check{\P}_{\SA\dot\x}$, and $\check{\P}_{\SA\dot\x,\SA\ddot\x}$.

The density $p(\SA\dot{\x},\SA\ddot{\x}|\ddot{\Y})$ is normally not known. However, from the definition of conditional probability
\begin{equation*}
p(\SA\dot\x,\SA\ddot\x|\ddot\Y)=p(\SA\dot\x|\SA\ddot\x,\ddot\Y)p(\SA\ddot\x|\ddot\Y).
\end{equation*}
$p(\SA\dot\x|\SA\ddot\x,\ddot\Y)$ can be obtained from the state equation via
\begin{equation*}
p(\SA\dot{\x}|\SA\ddot{\x})=\int\delta(\SA\dot\x-f^a(\SA\ddot\x,\dot\w))p(\dot\w)\d\dot\w.
\end{equation*}
Inserting into the expressions for the marginal moments and integrating out $\SA\dot\x$ give
\begin{gather*}
\SA\check{\dot\x}=\int{f}^a(\SA\ddot\x,\dot\w)p(\dot\w,\SA\ddot\x|\ddot\Y)\d([\SA\ddot{\x};\dot\w])\\
\check{\P}_{\SA\dot\x}=\int\big(f^a(\SA\ddot\x,\dot\w)-\SA\check{\dot\x}\big)^{\bot\top}p(\dot\w,\SA\ddot\x|\ddot\Y)\d([\SA\ddot{\x};\dot\w])\\
\check{\P}_{\SA\dot\x,\SA\ddot\x}=\int\big(f^a(\SA\ddot\x,\dot\w)-\SA\check{\dot\x}\big)(\SA\ddot\x-\SA\check{\ddot\x})^{\top}p(\dot\w,\SA\ddot\x|\ddot\Y)\d([\SA\ddot{\x};\dot\w])
\end{gather*}
where, for uniformity, we have written $p(\dot\w)p(\SA\ddot\x|\ddot\Y)=p(\dot\w,\SA\ddot\x|\ddot\Y)$. Consequently, the integration over $\SA\dot\x$ is replaced with the integration over $\dot\w$ and the density $p(\SA\dot{\x},\SA\ddot{\x}|\ddot{\Y})$ is replaced with the known density $p(\dot\w,\SA\ddot{\x}|\ddot{\Y})$.

\subsection{Marginalization for update}
Let $\TC\in\mathbb{R}^{n,b}$ and $\TD\in\mathbb{R}^{n,n-b}$ such that $\T^{\invsign}=[\TC\:\TD]$. Further, let $\N=\TB\check{\dot{\P}}\TA^\top(\TA\check{\dot{\P}}\TA^\top)^{\invsign}$. Note that $p(\mathbf{T}\dot{\x}|\dot{\Y})=p(\TB\dot{\x}|\TA\dot{\x},\!\dot{\Y})p(\TA\dot{\x}|\dot{\Y})$. Assumptions \emph{a} and $\mathcal{N}\!$, respectively, imply $p(\TB\dot{\x}|\TA\dot{\x},\dot{\Y})=p(\TB\dot{\x}|\TA\dot{\x},\ddot{\Y})=\mathcal{N}_{\TB\dot{\x}}\big(\TB\check{\dot{\x}}+\N(\TA\dot{\x}-\TA\check{\dot\x}),\TB\check{\dot{\P}}\TB^\top-\N\TA\check{\dot{\P}}\TB^\top\big)$. Finally, introduce the marginal moments
\begin{gather*}
\TA\hat{\dot\x}=\int\TA\dot{\x} p(\TA\dot{\x}|\dot{\Y})\d(\TA\dot{\x})\\
\hat{\P}_{\TA\dot{\x}}=\int(\TA\dot{\x}-\TA\hat{\dot\x})^{\bot\top}p(\TA\dot{\x}|\dot{\Y})\d(\TA\dot{\x}).
\end{gather*}
Together, this gives the updated mean
\begin{align*}
\hat{\dot\x}&=\int\dot\x p(\dot\x|\dot\Y)d\dot\x\\
             &=\T^{\invsign}\int\T\dot\x p(\T\dot\x|\dot\Y)\d(\T\dot\x)\\
             &=\TC\int\TA\dot\x p(\TA\dot\x|\dot\Y)\d(\TA\dot\x)
              +\TD\int\TB\dot\x p(\TB\dot\x|\TA\dot\x)p(\TA\dot\x|\dot\Y)\d(\T\dot\x)\\
             &=\TC\TA\hat{\dot\x}
              +\TD\int\left(\TB\check{\dot{\x}}+\N(\TA\dot{\x}-\TA\check{\dot\x})\right)p(\TA\dot\x|\dot\Y)\d(\TA\dot\x)\\
             &=\TC\TA\hat{\dot\x}+\TD\big(\TB\check{\dot\x}+\N(\TA\hat{\dot\x}-\TA\check{\dot\x})\big).
\end{align*}
Further, note that $\TB\hat{\dot\x}=\TB\TC\TA\hat{\dot\x}+\TB\TD\big(\TB\check{\dot\x}+\N(\TA\hat{\dot\x}-\TA\check{\dot\x})\big)=\TB\check{\dot\x}+\N(\TA\hat{\dot\x}-\TA\check{\dot\x})$. Together, this gives the updated covariance
\begin{align*}
\hat{\dot\P}&=\int(\dot\x-\hat{\dot\x})^{\bot\top}p(\dot\x|\dot\Y)\d\dot\x\\
            &=\int\T^{\invsign}\T(\dot\x-\hat{\dot\x})^{\bot\top}\T^\top\T^{\invtrans}p(\T\dot\x|\dot\Y)\d(\T\dot\x)\\
            &=\int\Big(\TC(\TA\dot\x-\TA\hat{\dot\x})^{\bot\top}\TC^\top+\TD(\TB\dot\x-\TB\hat{\dot\x})^{\top\bot}\SD^\top\\ &\quad\quad\quad\quad+(\TC\TA(\dot\x-\hat{\dot\x})^{\bot\top}\TB^\top\TD^\top)^{\bot+\top}\Big)p(\T\dot\x|\dot\Y)\d(\T\dot\x)\\
            &=\int\TC(\TA\dot\x-\TA\hat{\dot\x})^{\bot\top}\TC^\top{p}(\TA\dot\x|\dot\Y)\d(\TA\dot\x)\\
            &\:+\int\Bigg(\TD(\TB\dot\x\!-\!\TB\hat{\dot\x})^{\bot\top}\TD^\top\!\!+\!\Big(\TC(\TA\dot\x\!-\!\TA\hat{\dot\x})(\TB\dot\x\!-\!\TB\hat{\dot\x})^\top\TD^\top\Big)^{\bot+\top}\Bigg)\\
            &\:\:\:\mathcal{N}_{\TB\dot{\x}}\big(\TB\check{\dot{\x}}\!+\!\N(\TA\dot{\x}\!-\!\TA\check{\dot\x}),\TB\check{\dot{\P}}\TB^\top\!\!-\!\N\TA\check{\dot{\P}}\TB^\top\big)p(\TA\dot\x|\dot\Y)\d(\TA\dot\x)\d(\TB\dot\x)\\
            &=\TC\hat{\P}_{\TA\dot\x}\TC^\top\\
            &\!\!\!+\!\int\TD\Big(\TB\check{\dot{\P}}\TB^\top\!\!-\!\N\TA\check{\dot{\P}}\TB^\top+\N(\TA\dot\x-\TA\hat{\dot\x})^{\bot\top}\N^\top\Big)\TD^\top\!p(\TA\dot\x|\dot\Y)\d(\TA\dot\x)\\
            &\!\!\!+\!\!\int\!\Big(\TC(\TA\dot\x\!-\!\TA\hat{\dot\x})(\TA\check{\dot\x}\!+\!\N(\TA\dot\x\!-\!\TA\check{\dot\x})\!-\!\TA\hat{\dot\x})^\top\TD^\top\Big)^{\bot+\top}\!\!\!p(\TA\dot\x|\dot\Y)\d(\TA\dot\x)\\
            &\!\!\!\!\!=\!\TC\hat{\P}_{\TA\dot\x}\!\TC^\top\!\!\!+\!\TD(\TB\check{\dot{\P}}\TB^\top\!\!-\!\N\TA\check{\dot{\P}}\TB^\top\!\!\!+\!\N\hat{\P}_{\TA\dot\x}\N^\top)\TD^\top\!\!\!+\!(\TC\hat{\P}_{\TA\dot\x}\N^\top\!\TD^{\!\top})^{\bot\!+\!\top}\!\!.
\end{align*}
Consequently, with assumptions $\mathcal{N}$ and \emph{a}, obtaining $p(\dot\x|\dot\Y)$ from $p(\dot\x|\ddot\Y)$, i.e. the update, boils down to calculating the marginal moments $\TA\hat{\dot\x}$ and $\hat{\P}_{\TA\dot{\x}}$.

The density $p(\TA\dot{\x}|\dot{\Y})$ is not in general known. However, from Bayes' theorem and the Markovian assumption of the state space model
\begin{equation*}
p(\TA\dot{\x}|\dot{\Y})={\nu}p(\dot\y|\TA\dot\x)p(\TA\dot\x|\ddot\Y)
\end{equation*}
where $\nu$ is a normalization constant. $p(\TA\dot\x|\ddot\Y)$ is the known prior and $p(\dot\y|\TA\dot\x)$ is the likelihood function. If the likelihood function is known, up to scale, then the the marginal moments can directly be expressed as
\begin{gather*}
\TA\hat{\dot\x}=\nu\int\TA\dot\x{p}(\dot\y|\TA\dot\x)p(\TA\dot\x|\ddot\Y)\d(\TA\dot\x)\\
\hat{\P}_{\TA\dot\x}=\nu\int(\TA\dot\x-\TA\hat{\dot\x})^{\bot\top}p(\dot\y|\TA\dot\x)p(\TA\dot\x|\ddot\Y)\d(\TA\dot\x).
\end{gather*}
If the likelihood function is not known, it can be obtained from the output model $h^a(\cdot,\dot{\v})$ and the density $p(\dot{\v})$ via
\begin{equation*}
p(\dot{\y}|\TA\dot\x)=\int\delta(\dot\y-h^a(\TA\dot\x,\dot\v))p(\dot\v)\d\dot\v.
\end{equation*}
Inserting into the marginal moment integrals yields
\begin{gather*}
\TA\hat{\dot\x}=\nu\int\TA\dot\x\delta(\dot\y-h^a(\TA\dot\x,\dot\v))p(\dot\v)p(\TA\dot\x|\ddot\Y)\d(\TA\dot\x)\d\dot\v\\
\hat{\P}_{\TA\dot\x}=\nu\int(\TA\dot\x-\TA\hat{\dot\x})^{\bot\top}\delta(\dot\y-h^a(\TA\dot\x,\dot\v))p(\dot\v)p(\TA\dot\x|\ddot\Y)\d(\TA\dot\x)\d\dot\v.
\end{gather*}
Unfortunately, this increases the dimensionality of the integrals and adds the Dirac factor which complicates the evaluation. However, as will be shown in Section~\ref{sec:Marginal_moments_approximation}, the moments are still evaluable with standard statistical methods.

\section{Exploiting conditionally additive noise structure}\label{sec:exploiting_conditionally_additive_noise_structure}
If, in addition to the assumptions $\mathcal{N}$ and \emph{a}, assumption \emph{b} is made, the calculation of the marginal moments $\SA\check{\dot\x}$, $\check{\P}_{\SA\dot\x,\SA\ddot\x}$, $\check{\P}_{\SA\dot\x}$, $\TA\hat{\dot\x}$, and $\hat{\P}_{\TA\dot\x}$ can be facilitated by marginalizing out $\SA\dot\x$ (without introducing $\dot\w$) and an explicit likelihood function. For brevity, hereinafter let $\G^b=\G^b(\SA\ddot\x)$ and $\J^b=\J^b(\SA\ddot\x)$, while retaining implicit dependencies.

\subsection{Current state marginalization}
From the definition of conditional probability $p(\SA\dot\x,\SA\ddot\x|\ddot\Y)=p(\SA\dot\x|\SA\ddot\x,\ddot\Y)p(\SA\ddot\x|\ddot\Y)$. Now, given $\SA\ddot{\x}$, $\SA\dot{\x}$ is an affine transformation of $\dot\w$ and consequently
\begin{gather*}
\begin{aligned}
\SA\check{\dot\x}&=\int\SA\dot{\x}p(\SA\dot{\x}|\SA\ddot{\x},\ddot{\Y})p(\SA\ddot{\x}|\ddot{\Y})\d(\SA\ddot{\x})\d(\SA\dot{\x})\\
                 &=\int{f}^b(\SA\ddot{\x})p(\SA\ddot{\x}|\ddot{\Y})\d(\SA\ddot{\x})
\end{aligned}\\
\begin{aligned}
\check{\P}_{\SA\dot\x}&=\int(\SA\dot{\x}-\SA\check{\dot\x})^{\bot\top}p(\SA\dot{\x}|\SA\ddot{\x},\ddot{\Y})p(\SA\ddot{\x}|\ddot{\Y})\d(\SA\ddot{\x})\d(\SA\dot{\x})\\
                      &=\int\big(\G^b\P_{\dot\w}{\G^b}^\top+(f^b(\SA\ddot{\x})-\SA\check{\dot\x})^{\bot\top}\big)p(\SA\ddot{\x}|\ddot{\Y})\d(\SA\ddot{\x})
\end{aligned}\\
\begin{aligned}
\check{\P}_{\SA\dot\x,\SA\ddot\x}&=\int(\SA\dot\x-\SA\check{\dot\x})(\SA\ddot\x-\SA\check{\ddot\x})^\top{p}(\SA\dot{\x}|\SA\ddot{\x},\ddot{\Y})p(\SA\ddot{\x}|\ddot{\Y})\d(\SA\dot{\x})\\
                                  &=\int(f^b(\SA\ddot\x)-\SA\check{\dot\x})(\SA\ddot\x-\SA\check{\ddot\x})^\top{p}(\SA\ddot{\x}|\ddot{\Y})\d(\SA\ddot{\x}).
\end{aligned}
\end{gather*}
In contrast to the general case, $p(\dot\w)$ can be integrated out and therefore, the dimensionality of the integrals is reduced.

\subsection{Explicit likelihood function}
Again, from Bayes' theorem and the Markovian assumption of the state space model, $p(\TA\dot{\x}|\dot{\Y})=\nu{p}(\dot\y|\TA\dot\x)p(\TA\dot\x|\ddot\Y)$ where $\nu$ is a normalization constant. Given $\TA\dot\x$, $\dot\y$ is an affine transformation of $\dot\v$. Consequently, the likelihood function
\begin{equation*}
p(\dot\y|\TA\dot\x)=|\J^b|^{\invsign}p_{\dot\v}\big((\J^b)^{\invsign}(\dot\y-h^b(\TA\dot\x))\big)
\end{equation*}
where the subscript $p_{\dot\v}(\cdot)$ has been used to indicate that this is the pdf of $\dot{\v}$ and where $|\cdot|$ denotes the determinant.
Therefore
\begin{gather*}
\TA\hat{\dot\x}=\nu\int\TA\dot{\x}|\J^b|^{\invsign}p_{\dot\v}\big((\J^b)^{\invsign}(\dot\y-h^b(\TA\dot\x))\big)p(\TA\dot\x|\ddot\Y)\d(\TA\dot{\x})\\
\hat{\P}_{\TA\dot{\x}}\!=\!\nu\!\int\!(\TA\dot{\x}\!-\!\TA\hat{\dot\x})^{\bot\!\top}|\J^b|^{\invsign}p_{\dot\v}\big((\J^b)^{\invsign}(\dot\y\!-\!h^b(\TA\dot\x))\big)p(\TA\dot\x|\ddot\Y)\d(\TA\dot{\x}).
\end{gather*}
In contrast to the general case, $p(\dot\y|\TA\dot\x)$ is always known and the noise $\dot\v$ does not add to the dimensionality of the integrals. 
%

\section{Linear and Gaussian subspace marginalization}\label{sec:Linear_and_Gaussian_subspace_marginalization}
If, in addition to the assumptions $\mathcal{N}$ and \emph{a}, assumption \emph{c} is made, the calculation of the marginal moments $\SA\check{\dot\x}$, $\check{\P}_{\SA\dot\x,\SA\ddot\x}$, $\check{\P}_{\SA\dot\x}$, $\TA\hat{\dot\x}$, and $\hat{\P}_{\TA\dot\x}$ can be facilitated by further marginalizing out the linear subspaces. For brevity, hereinafter let $\F^c=\F^c(\SA^n\ddot\x)$, $\G^c=\G^c(\SA^n\ddot\x)$, $\H^c=\H^c(\TA^n\dot\x)$, and $\J^c=\J^c(\SA^n\ddot\x)$, while retaining implicit dependencies.

\subsection{Marginalization for prediction}
Let $\F^c=[\F^n;\F^l]$, $\G^c=[\G^n;\G^l]$, and $f^c(\cdot)=[f^n(\cdot);f^l(\cdot)]$ such that the divisions are consistent with the dimensions of $\SA=[\SA^n;\SA^l]$.
From the definition of conditional probability
\begin{align*}
p(\SA\dot\x,\SA\ddot\x|\ddot\Y)&\!=\!p(\SA^l\dot\x,\SA^l\ddot\x|\SA^n\dot\x,\SA^n\ddot\x,\ddot\Y)p(\SA^n\dot\x,\SA^n\ddot\x|\ddot\Y)\\
                               &\!=\!p(\SA^l\dot\x,\SA^l\ddot\x|\SA^n\dot\x,\SA^n\ddot\x,\ddot\Y)p(\SA^n\dot\x|\SA^n\ddot\x,\ddot\Y)p(\SA^n\ddot\x|\ddot\Y).
\end{align*}
The two first distributions can be derived in closed form.
Assumption $\mathcal{N}$ implies that
\begin{equation*}
p(\SA^l\ddot\x|\SA^n\ddot\x,\ddot\Y)=\mathcal{N}_{\SA^l\ddot\x}(\boldsymbol{\phi}(\SA^n\ddot\x),\boldsymbol{\Phi}) \end{equation*}
where
\begin{gather*} \boldsymbol{\phi}(\SA^n\ddot\x)=\SA^l\check{\ddot\x}+\SA^l\check{\ddot\P}{\SA^n}^\top(\SA^n\check{\ddot\P}{\SA^n}^\top)^{\invsign}(\SA^n\ddot\x-\SA^n\check{\ddot\x})\\
\boldsymbol{\Phi}=\SA^l\check{\ddot\P}{\SA^l}^\top-\SA^l\check{\ddot\P}{\SA^n}^\top(\SA^n\check{\ddot\P}{\SA^n}^\top)^{\invsign}\SA^n\check{\ddot\P}{\SA^l}^\top.
\end{gather*}
Together with assumptions $\mathcal{N}$ and \emph{c}, this implies
\begin{equation*}
p(\SA^n\dot\x|\SA^n\ddot\x,\ddot\Y)=\mathcal{N}_{\SA^n\dot\x}(\boldsymbol{\pi}(\SA^n\ddot\x),\boldsymbol{\Pi}(\SA^n\ddot\x))
\end{equation*}
where
\begin{align*}
\boldsymbol{\pi}(\SA^n\ddot\x)&=f^n(\SA^n\ddot\x)+\F^n\boldsymbol{\phi}(\SA^n\ddot\x)\\ \boldsymbol{\Pi}(\SA^n\ddot\x)&=\F^n\boldsymbol{\Phi}\F^n{}^\top+\G^n\P_{\dot\w}{\G^n}^\top.
\end{align*}
Further, the relations between $\SA^l\dot\x$, $\SA^l\ddot\x$, $\SA^n\dot\x$, and $\SA^n\ddot\x$ can be expressed as follows
\begin{align*}
\begin{bmatrix}
\SA^l\dot\x\\
\SA^l\ddot\x
\end{bmatrix}&=
\begin{bmatrix}
\F^l & \mathbf{0}\\
\mathbf{I} & \mathbf{0}
\end{bmatrix}
\begin{bmatrix}
\SA^l\ddot\x\\
\SA^l\dddot{\x}
\end{bmatrix}+
\begin{bmatrix}
f^l(\SA^n\ddot\x)\\
\mathbf{0}
\end{bmatrix}+
\begin{bmatrix}
\G^l\\
\mathbf{0}
\end{bmatrix}\dot\w\\
\SA^n\dot\x&=
\begin{bmatrix}
\mathbf{0}&\F^n
\end{bmatrix}
\begin{bmatrix}
\SA^l\dot\x\\
\SA^l\ddot\x
\end{bmatrix}+
f^n(\SA^n\ddot\x)+\G^n\dot\w.
\end{align*}
Given $\SA^n\dot\x$ and $\SA^n\ddot\x$, the former relation can be interpreted as a state equation and the latter as an output equation of a linear Gaussian system. $p(\SA^l\ddot\x|\SA^n\ddot\x,\ddot\Y)$ provides a Gaussian prior (note that $\SA^l\dddot\x$ is merely a placeholder since it does not effect the prediction) and therefore $p(\SA^l\dot\x,\SA^l\ddot\x|\SA^n\dot\x,\SA^n\ddot\x,\ddot\Y)$ is determined by the Kalman filter recursion, giving
\begin{gather*}
p(\SA^l\dot\x,\SA^l\ddot\x|\SA^n\dot\x,\SA^n\ddot\x,\ddot\Y)=\\
\begin{aligned}
\mathcal{N}_{\SA^l\dot\x,\SA^l\ddot\x}\Bigg(&
\begin{bmatrix}
f^l(\SA^n\ddot\x)+\F^l\boldsymbol{\phi}(\SA^n\ddot\x)+\boldsymbol{\Xi}^{ln}(\boldsymbol{\Xi}^n)^{\invsign}(\SA^n\dot\x-\boldsymbol{\pi}(\SA^n\ddot\x))\\
\boldsymbol{\phi}(\SA^n\ddot\x)+\boldsymbol{\Phi}{\F^n}^\top(\boldsymbol{\Xi}^n)^{\invsign}(\SA^n\dot\x-\boldsymbol{\pi}(\SA^n\ddot\x))
\end{bmatrix},\\&\:\:
\begin{bmatrix}
\boldsymbol{\Xi}^l-\boldsymbol{\Xi}^{ln}(\boldsymbol{\Xi}^n)^{\invsign}{\boldsymbol{\Xi}^{ln}}^\top & \F^l\boldsymbol{\Phi}-\boldsymbol{\Xi}^{ln}(\boldsymbol{\Xi}^n)^{\invsign}\F^n\boldsymbol{\Phi}\\
(\F^l\boldsymbol{\Phi}-\boldsymbol{\Xi}^{ln}(\boldsymbol{\Xi}^n)^{\invsign}\F^n\boldsymbol{\Phi})^\top & \boldsymbol{\Phi}-\boldsymbol{\Phi}{\F^n}^\top(\boldsymbol{\Xi}^n)^{\invsign}\F^n\boldsymbol{\Phi}
\end{bmatrix}
\Bigg)
\end{aligned}
\end{gather*}
where
\begin{align*}
\boldsymbol{\Xi}^l&=\F^l\boldsymbol{\Phi}{\F^l}^\top+\G^l\P_{\dot\w}{\G^l}^\top\\ \boldsymbol{\Xi}^{ln}&=\F^l\boldsymbol{\Phi}{\F^n}^\top+\G^l\P_{\dot\w}{\G^n}^\top\\ \boldsymbol{\Xi}^n&=\F^n\boldsymbol{\Phi}{\F^n}^\top+\G^n\P_{\dot\w}{\G^n}^\top.
\end{align*}
Consequently, $\SA^l\dot\x$ and $\SA^l\ddot\x$ and then $\SA^n\dot\x$ can be marginalized out of the marginal moments $\SA\check{\dot\x}$, $\check{\P}_{\SA\dot\x}$, and $\check{\P}_{\SA\dot\x,\SA\ddot\x}$.
Introduce
\begin{align*}
\boldsymbol{\alpha}(\SA^n\ddot\x)&=f^l(\SA^{n}\ddot\x)+\F^l\boldsymbol{\phi}(\SA^{n}\ddot\x)-\SA^l\check{\dot\x}\\ \boldsymbol{\beta}(\SA^n\ddot\x)&=\boldsymbol{\phi}(\SA^n\ddot\x)-\SA^l\check{\ddot\x}\\
\boldsymbol{\gamma}(\SA^n\dot\x)&=\boldsymbol{\pi}(\SA^n\ddot\x)-\SA^n\check{\dot\x}\\
\boldsymbol{\delta}(\SA^n\ddot\x)&=\SA^n\ddot\x-\SA^n\check{\ddot\x}\\
\P_{\boldsymbol{\alpha}}(\SA^n\ddot\x)&=\boldsymbol{\Xi}^l-\boldsymbol{\Xi}^{ln}(\boldsymbol{\Xi}^n)^{\invsign}{\boldsymbol{\Xi}^{ln}}^\top\\ \P_{\boldsymbol{\alpha},\boldsymbol{\beta}}(\SA^n\ddot\x)&=\F^l\boldsymbol{\Phi}-\boldsymbol{\Xi}^{ln}(\boldsymbol{\Xi}^n)^{\invsign}\F^n\boldsymbol{\Phi}\\
\boldsymbol{\epsilon}(\SA^n\ddot\x)&=f^l(\SA^{n}\ddot\x)+\F^l\boldsymbol{\phi}(\SA^{n}\ddot\x).
\end{align*}
Then, finally, straightforward calculations give
\begin{equation*}
\SA\check{\dot\x}=\int
\begin{bmatrix}
\boldsymbol{\pi}\\\boldsymbol{\epsilon}
\end{bmatrix}
p(\SA^n\ddot\x|\ddot\Y)\d(\SA^n\ddot\x)
\end{equation*}
\begin{equation*}
\:\:\check{\P}_{\SA\dot\x}=\int
\begin{bmatrix}
\boldsymbol{\Pi}+\boldsymbol{\gamma}\boldsymbol{\gamma}^\top&\boldsymbol{\gamma}\boldsymbol{\alpha}^\top+\boldsymbol{\Pi}\boldsymbol{\Xi}^\top\\
\boldsymbol{\alpha}\boldsymbol{\gamma}^\top+\boldsymbol{\Xi}\boldsymbol{\Pi} & \P_{\boldsymbol{\alpha}}+\boldsymbol{\alpha}\boldsymbol{\alpha}^\top+\boldsymbol{\Xi}\boldsymbol{\Pi}\boldsymbol{\Xi}^\top
\end{bmatrix}
p(\SA^n\ddot\x|\ddot\Y)\d(\SA^n\ddot\x)
\end{equation*}
\begin{equation*}
\check{\P}_{\SA\dot\x,\SA\ddot\x}=\int
\begin{bmatrix}
\boldsymbol{\gamma}\boldsymbol{\delta}^\top & \boldsymbol{\gamma}\boldsymbol{\beta}^\top+\boldsymbol{\Pi}\boldsymbol{\Gamma}^\top\\
\boldsymbol{\alpha}\boldsymbol{\delta}^\top & \P_{\boldsymbol{\alpha},\boldsymbol{\beta}}+\boldsymbol{\alpha}\boldsymbol{\beta}^\top+\boldsymbol{\Xi}\boldsymbol{\Pi}\boldsymbol{\Gamma}^\top
\end{bmatrix}
p(\SA^n\ddot\x|\ddot\Y)\d(\SA^n\ddot\x)
\end{equation*}
where $\boldsymbol{\Xi}=\boldsymbol{\Xi}^{ln}(\boldsymbol{\Xi}^n)^{\invsign}$ and $\boldsymbol{\Gamma}=\boldsymbol{\Phi}{\F^n}^\top(\boldsymbol{\Xi}^n)^{\invsign}$ and
where, for brevity, explicit dependencies on $\SA^n\ddot\x$ have been dropped. Consequently, the dimensionality of the moment integrals has been further reduced by the dimensionality of $\SA^l\ddot\x$.

\subsection{Marginalization for update}
From the definition of conditional probability
\begin{equation*}
p(\TA\dot\x|\dot\Y)=p(\TA^l\dot\x|\TA^n\dot\x,\dot\Y)p(\TA^n\dot\x|\dot\Y).
\end{equation*}
Assumption \emph{c} implies that, given $\TA^n\dot\x$, $\dot\y$ is linear in $\TA^l\dot\x$ with additive Gaussian noise. Let
\begin{gather*}
\boldsymbol{\psi}(\TA^n\dot\x)=\TA^l\check{\dot\x}+\TA^l\check{\dot\P}{\TA^n}^\top(\TA^n\check{\dot\P}{\TA^n}^\top)^{\invsign}(\TA^n\dot\x-\TA^n\check{\dot\x})\\ \boldsymbol{\Psi}=\TA^l\check{\dot\P}{\TA^l}^\top-\TA^l\check{\dot\P}{\TA^n}^\top(\TA^n\check{\dot\P}{\TA^n}^\top)^{\invsign}\TA^n\check{\dot\P}{\TA^l}^\top.
\end{gather*}
Then assumption $\mathcal{N}$ implies $p(\TA^l\dot\x|\TA^n\dot\x,\ddot\Y)\!=\!\mathcal{N}_{\TA^l\dot\x}\left(\boldsymbol{\psi}(\TA^n\dot\x), \boldsymbol{\Psi}\right)$.  Therefore, the conditioning on $\dot\y=\dot\Y\setminus\ddot\Y$ in $p(\TA^l\dot\x|\TA^n\dot\x,\dot\Y)$ is given by a Kalman filter update. Introduce
\begin{align*}
\varsigma(\TA^n\dot\x)&=h^c(\TA^n\dot\x)+\H^c\boldsymbol{\psi}(\TA^n\dot\x)\\ \Sigma(\TA^n\dot\x)&=\H^c\boldsymbol{\Psi}{\H^c}^\top+\J^c\P_{\dot\v}{\J^c}^\top\\
\boldsymbol{\omega}(\TA^n\dot\x)&=\boldsymbol{\psi}(\TA^n\dot\x)+\boldsymbol{\Psi}{\H^c}^\top\Sigma^{\invsign}(\TA^n\dot\x)\big(\dot\y-\varsigma(\TA^n\dot\x)\big)\\
\boldsymbol{\kappa}(\TA^n\dot\x)&=\boldsymbol{\omega}(\TA^n\dot\x)-\TA^l\hat{\dot\x}\\
\boldsymbol{\varkappa}(\TA^n\dot\x)&=\TA^n\dot\x-\TA^n\hat{\dot\x}\\
\boldsymbol{\Omega}(\TA^n\dot\x)&=(\I-\boldsymbol{\Psi}{\H^c}^\top\Sigma^{\invsign}(\TA^n\dot\x)\H^c)\boldsymbol{\Psi}.
\end{align*}
This gives
\begin{equation*}
p(\TA^l\dot\x|\TA^n\dot\x,\dot\Y)\!=\!\mathcal{N}_{\TA^l\dot\x}\left(\boldsymbol{\omega}(\TA^n\dot\x),\boldsymbol{\Omega}(\TA^n\dot\x)\right).
\end{equation*}
Consequently, $\TA^l\dot\x$ can be marginalized out of the marginal moments $\TA\hat{\dot\x}$ and $\hat{\P}_{\TA\dot{\x}}$. Again from Bayes' theorem $p(\TA^n\dot\x|\dot\Y)=\nu{p}(\dot\y|\TA^n\dot\x,\ddot\Y)p(\TA^n\dot\x|\ddot\Y)$, where $\nu$ is a normalization constant, and the likelihood function $p(\dot\y|\TA^n\dot\x,\ddot\Y)=\mathcal{N}_{\dot\y}(\varsigma(\TA^n\dot\x),\Sigma(\TA^n\dot\x))$. Together, this gives the marginal moments
\begin{equation*}
\TA\hat{\dot\x}=\nu\int
\begin{bmatrix}
\TA^n\dot{\x}\\
\boldsymbol{\omega}
\end{bmatrix}
\mathcal{N}_{\dot\y}(\varsigma,\Sigma)p(\TA^n\dot\x|\ddot\Y)\d(\TA^n\dot{\x})
\end{equation*}
\begin{equation*}
\hat{\P}_{\TA\dot{\x}}=\nu\int
\begin{bmatrix}
\boldsymbol{\varkappa}\boldsymbol{\varkappa}^\top&\boldsymbol{\varkappa}\boldsymbol{\kappa}^\top\\
\boldsymbol{\kappa}\boldsymbol{\varkappa}^\top&\boldsymbol{\Omega}+\boldsymbol{\kappa}\boldsymbol{\kappa}^\top
\end{bmatrix}
\mathcal{N}_{\dot\y}(\varsigma,\Sigma)p(\TA^n\dot\x|\ddot\Y)\d(\TA^n\dot{\x})
\end{equation*}
where, for brevity, explicit dependencies on $\TA^n\dot\x$ have been dropped. Consequently, the dimensionality of the moment integrals has been further reduced by that of $\TA^l\dot\x$.

\section{Closed form solutions for the linear case}\label{sec:Closed_form_solutions_for_the_linear_case}
If, in addition to $\mathcal{N}$ and \emph{a}, assumption \emph{d} is made, closed form solutions for the marginal moments can be derived. Substituting the formulas below back into the expressions for the moments gives the Kalman filter predictions and updates.

\subsection{Closed form solutions for prediction}
Substituting $\f^d+\F^d\SA\ddot\x+\G^d\dot\w$ for $f^a(\SA\ddot\x,\dot\w)$ in the formulas for the marginal moments given assumptions $\mathcal{N}$ and \emph{a}, gives
\begin{gather*}
\begin{aligned}
\SA\check{\dot\x}&=\int(\f+\F^d\SA\ddot\x+\G^d\dot\w)p(\dot\w,\SA\ddot\x|\ddot\Y)\d([\SA\ddot{\x};\dot\w])\\
&=\f^d+\F^d\SA\check{\ddot\x}
\end{aligned}\\
\begin{aligned}
\check{\P}_{\SA\dot\x}&=\int(\f+\F^d\SA\ddot\x+\G^d\dot\w-\SA\check{\dot\x})^{\bot\top}p(\dot\w,\SA\ddot\x|\ddot\Y)\d([\SA\ddot{\x};\dot\w])\\
&=\F^d\SA\check{\ddot\P}\SA^\top{\F^d}^\top+\G^d\P_{\dot\w}{\G^d}^\top
\end{aligned}\\
\begin{aligned}
\check{\P}_{\SA\dot\x,\SA\ddot\x}\!&=\!\!\!\!\:\int\!(\f\!\!\!\;+\!\F^d\!\!\:\SA\ddot\x\!+\!\G^d\dot\w\!-\!\!\!\;\SA\check{\dot\x})(\SA\ddot\x\!-\!\SA\check{\ddot\x})^{\top}\!\!\!\;p(\dot\w,\!\!\:\SA\ddot\x|\ddot\Y)\d([\!\!\;\SA\ddot{\x};\!\!\:\dot\w\!\!\;])\\
&=\F^d\SA\check{\ddot\P}\SA^\top.
\end{aligned}
\end{gather*}

\subsection{Closed form solutions for update}
With linear and Gaussian output, the marginal moments are given by the Kalman filter update formulas
\begin{align*}
\TA\hat{\dot\x}\!\!\:&=\!\!\:\TA\check{\dot\x}\!+\!\TA\check{\dot\P}\TA^\top{\H^d}^\top\!(\H^d\TA\check{\dot\P}\TA^\top{\H^d}^\top\!\!\!\!\:+\!\J^d\P_{\dot\v}{\J^d}^\top\!\!\:)^{\invsign} \!\!\;(\dot\y\!-\!\h^d\!\!\!\:-\!\H^d\!\!\;\TA\check{\dot\x})\\
\hat\P_{\TA\dot\x}\!\!\:&=\!\!\:(\I\,-\TA\check{\dot\P}\TA^\top{\H^d}^\top\!(\H^d\TA\check{\dot\P}\TA^\top{\H^d}^\top\!\!\!\!\:+\!\J^d\P_{\dot\v}{\J^d}^\top\!\!\:)^{\invsign}\H^d)\TA\check{\dot\P}\TA^\top.
\end{align*}

\section{Numerical moment integral approximations}~\label{sec:Marginal_moments_approximation}
With assumption $\mathcal{N}$, the Bayesian filtering is reduced to approximating the moments $\{\check{\dot\x},\check{\dot\P},\hat{\dot\x},\hat{\dot\P}\}$. With  assumption \emph{a},
\begin{align*}
\check{\dot{\x}}&=\SC\SA\check{\dot\x}+\SD\SB\check{\ddot{\x}}\\
\check{\dot{\P}}&=\SC\check{\P}_{\SA\dot\x}\SC^\top+\SD\SB\check{\ddot\P}\SB^\top\SD^\top+(\SC\check{\P}_{\SA\dot\x,\SA\ddot\x}\M^\top\SD^\top)^{\bot+\top}\\
\hat{\dot\x}&=\TC\TA\hat{\dot\x}+\TD\big(\TB\check{\dot\x}+\N(\TA\hat{\dot\x}-\TA\check{\dot\x})\big)\\
\hat{\dot\P}&=\!\TC\hat{\P}_{\TA\dot\x}\!\TC^{\!\top}\!\!\!+\!\TD(\TB\check{\dot{\P}}\TB^{\!\top}\!\!\!-\!\N\TA\check{\dot{\P}}\TB^{\!\top}\!\!\!+\!\N\hat{\P}_{\TA\dot\x}\N^{\!\top})\TD^{\!\top}\!\!\!+\!(\TC\hat{\P}_{\TA\dot\x}\N^{\!\top}\!\TD^{\!\top})^{\bot\!+\!\top}
\end{align*}
and the approximations can be focused to the marginal moments $\{\SA\check{\dot\x},\check{\P}_{\SA\dot\x},\check{\P}_{\SA\dot\x,\SA\ddot\x},\TA\hat{\dot\x},\hat{\P}_{\TA\dot{\x}}\}$. %
Assumptions \emph{b} and \emph{c} enable us to further narrow the focus, together giving capability 1.

All non-closed-form formulas derived for the marginal moments are on the form
\begin{equation*}
\m=\int{d}(\cdot)p(\cdot|\ddot\Y)\d(\cdot)
\end{equation*}
where $d(\cdot)$ is a function of the subspace indicated by the marginal density $p(\cdot|\ddot\Y)$ and $\d(\cdot)$ is a differential of the corresponding subspace.
Given assumptions $\mathcal{N}$ and \emph{a}
\begin{gather*}
p(\cdot|\ddot\Y)=p(\dot\w,\SA\ddot\x|\ddot\Y)=p(\dot\w)\mathcal{N}_{\SA\x}(\SA\check{\ddot{\x}},\SA\check{\ddot{\P}}\SA^\top)\quad\text{and}\\
\begin{aligned}
d(\dot\w,\SA\ddot\x)&=f^a(\SA\ddot\x,\dot\w)                                &&\Rightarrow\m=\SA\check{\dot\x}\\
d(\dot\w,\SA\ddot\x)&=(f^a(\SA\ddot\x,\dot\w)-\SA\check{\dot\x})^{\bot\top} &&\Rightarrow\m=\check{\P}_{\SA\dot\x}\\
d(\dot\w,\SA\ddot\x)&=(f^a(\SA\ddot\x,\dot\w)-\SA\check{\dot\x})(\SA\ddot\x-\SA\check{\ddot\x})^{\top} &&\Rightarrow\m=\check{\P}_{\SA\dot\x,\SA\ddot\x}.
\end{aligned}
\end{gather*}
and additionally given $p(\dot\y|\TA\dot\x)$
\begin{gather*}
p(\cdot|\ddot\Y)=\mathcal{N}_{\TA\dot{\x}}(\TA\check{\dot\x},\TA\check{\dot{\P}}\TA^\top)\quad\text{and}\\
\begin{aligned}
d(\TA\dot\x)&=\nu\TA\dot{\x}p(\dot\y|\TA\dot\x)
&&\Rightarrow\m=\TA\hat{\dot\x}\\
d(\TA\dot\x)&=\nu(\TA\dot{\x}\!-\!\TA\hat{\dot\x})^{\bot\!\top}p(\dot\y|\TA\dot\x) &&\Rightarrow\m=\hat{\P}_{\TA\dot{\x}}
\end{aligned}
\end{gather*}
otherwise (if $p(\dot\y|\TA\dot\x)$ is not available)
\begin{gather*}
\begin{aligned}
d(\TA\dot\x)&=\nu\TA\dot{\x}\left({\textstyle\int}\delta(\dot\y-h^a(\TA\dot\x,\dot\v))p(\dot\v)\d\dot\v\right)
&&\Rightarrow\m=\TA\hat{\dot\x}\\
d(\TA\dot\x)&=\nu(\TA\dot{\x}\!-\!\TA\hat{\dot\x})^{\bot\!\top}\!\left({\textstyle\int}\delta(\dot\y\!-\!h^a(\TA\dot\x,\!\dot\v))p(\dot\v)\d\dot\v\right)\!\! &&\Rightarrow\m=\hat{\P}_{\TA\dot{\x}}.
\end{aligned}
\end{gather*}
Given assumptions $\mathcal{N}$, \emph{a}, and \emph{b}
\begin{gather*}
p(\cdot|\ddot\Y)=\mathcal{N}_{\SA\x}(\SA\check{\ddot{\x}},\SA\check{\ddot{\P}}\SA^\top)\quad\text{and}\\
\begin{aligned}
d(\SA\ddot\x)&=f^b(\SA\ddot{\x})
&&\Rightarrow\m=\SA\check{\dot\x}\\
d(\SA\ddot\x)&=\G^b\P_{\dot\w}{\G^b}^\top+(f^b(\SA\ddot{\x})-\SA\check{\dot\x})^{\bot\top} &&\Rightarrow\m=\check{\P}_{\SA\dot\x}\\
d(\SA\ddot\x)&=(f^b(\SA\ddot\x)-\SA\check{\dot\x})(\SA\ddot\x-\SA\check{\ddot\x})^\top &&\Rightarrow\m=\check{\P}_{\SA\dot\x,\SA\ddot\x}
\end{aligned}
\end{gather*}
and
\begin{gather*}
p(\cdot|\ddot\Y)=\mathcal{N}_{\TA\dot{\x}}(\TA\check{\dot\x},\TA\check{\dot{\P}}\TA^\top)\quad\text{and}\\
\begin{aligned}
d(\TA\dot\x)&=\!\nu\TA\dot{\x}|\J^b|^{\invsign}p_{\dot\v}\big((\J^b)^{\invsign}(\dot\y-h^b(\TA\dot\x))\big)\!\!\!
&&\Rightarrow\!\m=\TA\hat{\dot\x}\\
d(\TA\dot\x)&=\!\nu(\TA\dot{\x}\!-\!\TA\hat{\dot\x})^{\bot\!\top}|\J^b|^{\invsign}p_{\dot\v}\big((\J^b)^{\invsign}(\dot\y\!-\!h^b(\TA\dot\x))\big)\!\!\! &&\Rightarrow\!\m=\hat{\P}_{\TA\dot{\x}}\!\!\;.
\end{aligned}
\end{gather*}
Finally, given assumptions $\mathcal{N}$, \emph{a}, and \emph{c}
\begin{gather*}
p(\cdot|\ddot\Y)=\mathcal{N}_{\SA^n\x}(\SA^n\check{\ddot{\x}},\SA^n\check{\ddot{\P}}{\SA^n}^\top)\quad\text{and}\\
\begin{aligned}
d(\SA\ddot\x)&=[\boldsymbol{\pi};\boldsymbol{\epsilon}]
&&\Rightarrow\m=\SA\check{\dot\x}\\
d(\SA\ddot\x)&=
\begin{bmatrix}
\boldsymbol{\Pi}+\boldsymbol{\gamma}\boldsymbol{\gamma}^\top&\boldsymbol{\gamma}\boldsymbol{\alpha}^\top+\boldsymbol{\Pi}\boldsymbol{\Xi}^\top\\
\boldsymbol{\alpha}\boldsymbol{\gamma}^\top+\boldsymbol{\Xi}\boldsymbol{\Pi} & \P_{\boldsymbol{\alpha}}+\boldsymbol{\alpha}\boldsymbol{\alpha}^\top+\boldsymbol{\Xi}\boldsymbol{\Pi}\boldsymbol{\Xi}^\top
\end{bmatrix}
&&\Rightarrow\m=\check{\P}_{\SA\dot\x}\\
d(\SA\ddot\x)&=
\begin{bmatrix}
\boldsymbol{\gamma}\boldsymbol{\delta}^\top & \boldsymbol{\gamma}\boldsymbol{\beta}^\top+\boldsymbol{\Pi}\boldsymbol{\Gamma}^\top\\
\boldsymbol{\alpha}\boldsymbol{\delta}^\top & \P_{\boldsymbol{\alpha},\boldsymbol{\beta}}+\boldsymbol{\alpha}\boldsymbol{\beta}^\top+\boldsymbol{\Xi}\boldsymbol{\Pi}\boldsymbol{\Gamma}^\top
\end{bmatrix}
&&\Rightarrow\m=\check{\P}_{\SA\dot\x,\SA\ddot\x}
\end{aligned}
\end{gather*}
and
\begin{gather*}
p(\cdot|\ddot\Y)=\mathcal{N}_{\TA^n\dot{\x}}(\TA^n\check{\dot{\x}},\TA^n\check{\dot{\P}}{\TA^n}^\top)\quad\text{and}\\
\begin{aligned}
d(\TA\dot\x)&=\nu\,\big[\TA^n\dot{\x};\boldsymbol{\omega}\big]\,\mathcal{N}_{\dot\y}(\varsigma,\Sigma)
&&\Rightarrow\m=\TA\hat{\dot\x}\\
d(\TA\dot\x)&=\nu
\begin{bmatrix}
\boldsymbol{\varkappa}\boldsymbol{\varkappa}^\top&\boldsymbol{\varkappa}\boldsymbol{\kappa}^\top\\
\boldsymbol{\kappa}\boldsymbol{\varkappa}^\top&\boldsymbol{\Omega}+\boldsymbol{\kappa}\boldsymbol{\kappa}^\top
\end{bmatrix}
\mathcal{N}_{\dot\y}(\varsigma,\Sigma) &&\Rightarrow\m=\hat{\P}_{\TA\dot{\x}}.
\end{aligned}
\end{gather*}
Suitable numerical integration methods for $\m$ come in a range of different flavors such as Gauss-Hermit quadrature, importance sampling and Monte Carlo methods. All methods may be viewed as a set of samples $\r^{(i)}$  and weights $w^{(i)}$ giving
\begin{equation*}
\m\approx\bigg(\medop\sum_i{w}^{(i)}\bigg)^{\invsign}\medop\sum_id(\r^{(i)})w^{(i)}.
\end{equation*}
The suitable $\r^{(i)}$ depends on $d(\cdot)$ and the method of choice.

Given only assumptions $\mathcal{N}$ and \emph{a}, the above formula cannot be directly applied to the update since the integral within $d(\cdot)$ cannot typically be evaluated analytically and in general we have to resort to density estimation. With a kernel $\mathcal{K}_{\mathbf{Z}(\cdot)}(\cdot)$, a bandwidth matrix $\mathbf{Z}(\cdot)$, samples $\v^{(i,j)}$, and weights ${v}^{(i,j)}$,
\begin{align*}
\int\delta(\dot\y-h^a(\r^{(i)},\dot\v))&p(\dot\v)\d\dot\v\\ \approx\bigg(&\medop\sum_j{v}^{(i,j)}\bigg)^{\invsign}\medop\sum_j\mathcal{K}_{\mathbf{Z}(\r^{(i)})}\big(\dot{\y}-h^a(\s^{(i)},\v^{(i,j)})\big){v}^{(i,j)}.
\end{align*}
Substituting this back into $d(\cdot)$ enables use of numerical integration as described above. Again, the suitable $\v^{(i,j)}$ and $\mathcal{K}_{\mathbf{Z}(\cdot)}(\cdot)$ depends on $h^a(\cdot,\cdot)$ and $p(\dot\v)$ and the method of choice.

An alternative to the density estimation, for assumptions $\mathcal{N}$ and \emph{a}, is a parametric Kalman update. Extending assumption $\mathcal{N}$ to a jointly Gaussian assumption
\begin{equation*}
p(\TA\dot\x,\dot\y|\ddot\Y)=\mathcal{N}_{\TA\dot\x,\dot\y}\left(
\begin{bmatrix}
\check{\TA\dot\x}\\\check{\dot\y}
\end{bmatrix},
\begin{bmatrix}
\check{\P}_{\TA\dot\x} & \check{\P}_{\TA\dot\x,\dot\y}\\
\check{\P}_{\TA\dot\x,\dot\y}^\top & \check{\P}_{\dot\y}
\end{bmatrix}
\right)
\end{equation*}
implies
\begin{align*}
\TA\hat{\dot\x}&=\TA\check{\dot\x}+\check{\P}_{\TA\dot\x,\dot\y}(\check{\P}_{\dot\y})^{\invsign}\big(\dot\y-\check{\dot\y}\big)\\
\hat{\dot\P}_{\TA\dot\x}&=\check{\dot\P}_{\TA\dot\x}-\check{\P}_{\TA\dot\x,\dot\y}(\check{\P}_{\dot\y})^{\invsign}\check{\P}_{\TA\dot\x,\dot\y}^\top
\end{align*}
which are the well-known Kalman filter update formulas. Consequently, the problem is changed to that of approximating
\begin{gather*}
\check{\dot\y}=\int h^a(\TA\dot\x,\dot\v)p(\dot\v)p(\TA\dot\x|\ddot\Y)\d(\TA\dot\x)\d\dot\v\\
\check{\P}_{\dot\y}=\int(\check{\dot\y}-h^a(\TA\dot\x,\dot\v))^{\bot\top}p(\dot\v)p(\TA\dot\x|\ddot\Y)\d(\TA\dot\x)\d\dot\v\\
\check{\P}_{\TA\dot\x,\dot\y}=\int(\TA\dot\x-\TA\check{\dot\x})(\check{\dot\y}-h^a(\TA\dot\x,\dot\v))p(\dot\v)p(\TA\dot\x|\ddot\Y)\d(\TA\dot\x)\d\dot\v,
\end{gather*}
which are amenable to numerical integration as described above. Hence, the density estimation is avoided at the cost of the extended Gaussian assumption, giving capability 1 to nonlinear Kalman filters applied to the setup of assumption \emph{a}.

\section{Introducing linearization and Gaussian approximations}
To focus the computational effort on the significantly nonlinear, non-additive, and non-Gaussian model parts one may exploit linearization and Gaussian approximation of benign parts of the models. For simplicity, assume zero mean noise (and hence a zero noise linearization point) and restate assumptions \emph{b-c}:

\vspace{3mm}
\noindent\emph{b') Conditionally additive noise:}
Assumption \emph{b} may alternatively be expressed as
\begin{align*}
f^a(\SA\ddot\x,\dot\w)&={f}^a(\SA\ddot\x,\O)+\left.\frac{\partial{f}^a}{\partial\dot\w}\right|_{(\SA\ddot\x,\O)}\dot\w\\
h^a(\TA\dot\x,\dot\v)&={h}^a(\TA\dot\x,\O)+\left.\frac{\partial{h}^a}{\partial\dot\v}\right|_{(\TA\dot\x,\O)}\dot\v,
\end{align*}
giving
\begin{gather*}
f^b(\SA\ddot\x)=f^a(\SA\ddot\x,\O),\:\: h^b(\TA\dot\x)=h^a(\TA\dot\x,\O)\\ \G^b=\left.\frac{\partial{f}^a}{\partial\dot\w}\right|_{(\SA\ddot\x,\O)}\quad\text{and}\quad   \J^b=\left.\frac{\partial{h}^a}{\partial\dot\v}\right|_{(\TA\dot\x,\O)}.
\end{gather*}

\vspace{3mm}
\noindent\emph{c') Conditionally linear and Gaussian subspaces:}
Assumption \emph{c} may alternatively be expressed as
\begin{align*}
f^a(\SA\ddot\x,\!\dot\w)&\!=\!\!f^a\!\left(\left[\substack{\SA^n\ddot\x\\\SA^l\bar{\ddot\x}}\right]\!,\!\O\right)\!+\!\!\left.\frac{\partial{f}^a}{\partial\SA^l\ddot\x}\right|_{\left(\left[\substack{\SA^n\ddot\x\\\SA^l\bar{\ddot\x}}\right],\O\right)}\!\!(\SA^l\ddot\x\!-\!\SA^l\bar{\ddot\x})\!+\!\!\left.\frac{\partial{f}^a}{\partial\dot\w}\right|_{\left(\left[\substack{\SA^n\ddot\x\\\SA^l\bar{\ddot\x}}\right],\O\right)}\!\!\dot\w'\\
h^a(\TA\dot\x,\!\dot\v)&\!=\!h^a\!\left(\left[\substack{\TA^n\dot\x\\\TA^l\bar{\dot\x}}\right]\!,\!\O\right)\!+\!\!\left.\frac{\partial{h}^a}{\partial\TA^l\dot\x}\right|_{\left(\left[\substack{\TA^n\dot\x\\\TA^l\bar{\dot\x}}\right],\O\right)}\!\!(\TA^l\dot\x\!-\!\TA^l\bar{\dot\x})\!+\!\!\left.\frac{\partial{h}^a}{\partial\dot\v}\right|_{\left(\left[\substack{\TA^n\dot\x\\\TA^l\bar{\dot\x}}\right],\O\right)}\!\!\dot\v'\!,
\end{align*}
where $\dot\v'\sim\mathcal{N}_{\dot\v'}(\O,\P_{\dot\v})$ and $\dot\w'\sim\mathcal{N}_{\dot\w'}(\O,\P_{\dot\w})$ and $\SA^l\bar{\ddot\x}$ and $\TA^l\bar{\dot\x}$ are linearization points, giving
\begin{gather*}
\F^c\!=\!\left.\frac{\partial{f}^a}{\partial\SA^l\ddot\x}\right|_{\left(\left[\substack{\SA^n\ddot\x\\\SA^l\bar{\ddot\x}}\right],\O\right)},\; \G^c\!=\!\left.\frac{\partial{f}^a}{\partial\dot\w}\right|_{\left(\left[\substack{\SA^n\ddot\x\\\SA^l\bar{\ddot\x}}\right],\O\right)},\;
\H^c\!=\!\left.\frac{\partial{h}^a}{\partial\TA^l\dot\x}\right|_{\left(\left[\substack{\TA^n\dot\x\\\TA^l\bar{\dot\x}}\right],\O\right)},\\ \J^c\!=\!\left.\frac{\partial{h}^a}{\partial\dot\v}\right|_{\left(\left[\substack{\TA^n\dot\x\\\TA^l\bar{\dot\x}}\right],\O\right)},\quad\text{and}\quad
\begin{matrix}
f^c(\SA^n\ddot\x)\!=\!f^a\!\left(\left[\substack{\SA^n\ddot\x\\\SA^l\bar{\ddot\x}}\right]\!,\!\O\right)-\F^c\SA^l\bar{\ddot\x}\\ h^c(\TA^n\dot\x)=\!h^a\!\left(\left[\substack{\TA^n\dot\x\\\TA^l\bar{\dot\x}}\right]\!,\!\O\right)-\H^c\TA^l\bar{\dot\x}.
\end{matrix}
\end{gather*}

\vspace{3mm}
\noindent\emph{d') Linear and Gaussian subspaces:}
Assumption \emph{c} may alternatively be expressed as
\begin{align*}
f^a(\SA\ddot\x,\dot\w)&={f}^a(\SA\bar{\ddot\x},\O)+\left.\frac{\partial{f}^a}{\partial\SA\ddot\x}\right|_{(\SA\bar{\ddot\x},\O)}(\SA\ddot\x-\SA\bar{\ddot\x})+\left.\frac{\partial{f}^a}{\partial\dot\w}\right|_{(\SA\bar{\ddot\x},\O)}\dot\w\\
h^a(\TA\dot\x,\dot\v)&={h}^a(\TA\bar{\dot\x},\O)+\left.\frac{\partial{h}^a}{\partial\TA\dot\x}\right|_{(\TA\bar{\dot\x},\O)}(\TA\dot\x-\TA\bar{\dot\x})+\left.\frac{\partial{h}^a}{\partial\dot\v}\right|_{(\TA\bar{\dot\x},\O)}\dot\v',
\end{align*}
where $\dot\v'\sim\mathcal{N}_{\dot\v'}(\O,\P_{\dot\v})$ and $\SA\bar{\ddot\x}$ and $\TA\bar{\dot\x}$ are linearization points, giving
\begin{gather*}
\F^d\!\!\!\;=\!\!\!\:\left.\frac{\partial{f}^a}{\partial\SA\ddot\x}\right|_{(\SA\bar{\ddot\x},\O)}\!\!\!\;,\: \G^d\!\!\!\;=\!\!\!\:\left.\frac{\partial{f}^a}{\partial\dot\w}\right|_{(\SA\bar{\ddot\x},\O)}\!\!\!\;,\:
\H^d\!\!\!\;=\!\!\!\:\left.\frac{\partial{h}^a}{\partial\TA\dot\x}\right|_{(\TA\bar{\dot\x},\O)}\!\!\!\;,\: \J^d\!\!\!\;=\!\!\!\:\left.\frac{\partial{h}^a}{\partial\dot\v}\right|_{(\TA\bar{\dot\x},\O)}\\
\f^d={f}^a(\SA\bar{\ddot\x},\O)-\F^d\SA\bar{\ddot\x},\quad\text{and}\quad \h^d={h}^a(\TA\bar{\dot\x},\O)-\H^d\TA\bar{\dot\x}.
\end{gather*}

\noindent In principle, the assumptions $b/b'$, $c/c'$, and $d/d'$ are equivalent. In practice, the assumption versions \emph{b'}, \emph{c'} and \emph{d'} can be used to introduce model approximations by neglecting the first pair of relations in each assumption. This trades model approximations for lower dimensional moment integrals, and gives capability 2. Defining the above derivatives and covariances for state space models also makes it possible to automatically switch between different approximation modes.

The linearization points $\SA^l\bar{\ddot\x}$, $\TA^l\bar{\dot\x}$, $\SA\bar{\ddot\x}$, and $\TA\bar{\dot\x}$ are arbitrary but in practice they will normally be the marginal prior means $\SA^l\check{\ddot\x}$, $\TA^l\check{\dot\x}$, $\SA\check{\ddot\x}$, and $\TA\check{\dot\x}$.

\section{Decoupling modelling and system composition}\label{sec:system_composition}
A system model is composed of submodels. At least, as of the state space description, it is composed of a state model and an output model, but often there are multiple submodels of different states and outputs. Assume submodels
\begin{equation*}
\dot\x^i=f^i(\ddot\x^i,\dot\w^i),\quad\dot\y^j=h^j(\dot\x^j,\dot\v^j):\quad i,j\in\{1,2,\dots\}
\end{equation*}
where $\x^i$ and $\x^j$ are states of the properties determining the subsystem dynamics and the subsystem outputs, respectively. The system state $\x$ spans the properties of $\{\x^i,\x^j\}$. The properties of $\{\x^i,\x^j\}$ must overlap but individual model states will often not coincide with $\x$ and different models may use different ordering and units of the states.

The discrepancies between $\x$ and $\{\x^i,\x^j\}$ mean that the submodels cannot directly be compounded. Naively, this is solved by modifying the models, or using wrapper models, such that they are expressed in $\x$. However, this creates an unfortunate coupling between modelling and system composition, a maintenance nightmare, and prevents exploitation of subsystem structure (capability 1 and 2) since potential structure is hidden. Instead, using assumption \emph{a}, the composition can be performed by constructing $\SA$:s and $\TA$:s such that
\begin{equation*}
\x^i=\SA\x\quad\text{and}\quad\x^j=\TA\x.
\end{equation*}
The $\SA$:s and $\TA$:s will include scaling (unit transformations) and state selection but may also include further linear transformations~\cite{Nilsson2013}. This effectively decouples the modelling and the system composition giving capability 3.


\section{Relations to other filtering techniques and prior work}\label{sec:connections}
In the following subsection we describe how the presented results relate to the Kalman filters and the marginalized particle filters, and a few other directly related techniques. 

\subsection{Kalman filters}
The Kalman filters handle a setup with assumption $\mathcal{N}$ and the additional jointly Gaussian assumption
\begin{equation*}
p(\dot\x,\dot\y|\ddot\Y)=\mathcal{N}_{\dot\x,\dot\y}\left(
\begin{bmatrix}
\check{\dot\x}\\\check{\dot\y}
\end{bmatrix},
\begin{bmatrix}
\check{\dot\P} & \check{\P}_{\dot\x,\dot\y}\\
\check{\P}_{\dot\x,\dot\y}^\top & \check{\P}_{\dot\y}
\end{bmatrix}
\right).
\end{equation*}
(Note that Kalman filters can be formulated without the Gaussian assumption but the differences in practise are small.) The (Kalman filter) update rules for this case, in the active subspace, were presented in Section~\ref{sec:Marginal_moments_approximation}. Hence, the inactive subspace marginalization is directly applicable to nonlinear Kalman filters and they may be viewed as parametric methods for filtering given assumption $\mathcal{N}$, and in many cases additional assumptions such as additive noise.
See e.g.~\cite{Julier2000,Ito2000,Arasaratnam2009,Steinbring2016} in which different numerical approximation techniques for the moment integrals, together with additional model assumptions, are presented, giving the so called UKF, QKF, CKF and S${}^2$KF, respectively.

In the special case of assumptions $\mathcal{N}$ and \emph{d}, the jointly Gaussian condition is implied. Therefore, 
the formulas for this case give equivalent results to the Kalman filter. Further, a complete model linearization and Gaussian approximation, by using the latter part of assumption \emph{d'}, 
gives equivalent results to an explicitly linearized Kalman filter, e.g. the EKF.

Marginalization techniques for special cases of assumptions \emph{b} and \emph{c} have been presented in conjuncture with Kalman filters. In~\cite{Closas2010}, a marginalized QKF is presented in which a quadrature rule is used to compute predicted and updated mean and covariance for a nonlinear subspace. Subsequently, the mean, covariance and cross-covariance of the conditionally linear and Gaussian subspace are predicted and updated. In contrast, I present techniques for jointly calculating the means and covariances of the corresponding subspaces. Further, in~\cite{Briers2003,Chang2014} various marginalizations combined with a UKF are presented. However, the marginalizations are tailored to the jointly Gaussian Kalman filter setup. In~\cite{Huang2015}, marginalization of an inactive augmented noise state in the prediction is presented. In~\cite{Liu2011} a combination of linear approximations and numerical approximations for a UKF is presented. Finally, marginalization-like techniques for linear and Gaussian update subspaces in Kalman filters are found in the PUKF~\cite{Raitoharju2016}. Essentially, all the previous marginalization work rely on the jointly Gaussian assumption, and marginalization exploiting probably the most important structure, active subspaces, in the sense leading to the largest dimensionality reduction, is missing.

\subsection{Marginalized particle filters}
Marginalized, a.k.a. Rao-Blackwellized, particle filters, see e.g.~\cite{Schon2005}, are different in nature since they do not make assumption $\mathcal{N}$. However, they can only handle special cases of assumptions \emph{a}-\emph{c}. In general, marginalized particle filters require a \emph{static} conditional linear and Gaussian subspace,
\begin{equation*}
\SA^n_k=\TA^n_k=\SA^n=\TA^n
\end{equation*}
i.e. the marginalized particle filter handle a special case of assumption \emph{c}. (Throughout the results we have dropped indices $k$ of the subspace matrices but it has been implicit that they may vary with $k$.) The static subspace implies that
\begin{gather*}
p(\SA^l\dot\x|\SA^n\dot\x,\dot\Y)=\mathcal{N}_{\SA^l\dot\x}(\SA^l\hat{\dot\x},\hat{\dot\P})\\
p(\SA^l\dot\x|\SA^n\ddot\x,\dot\Y)=\mathcal{N}_{\SA^l\dot\x}(\SA^l\check{\dot\x},\check{\dot\P})
\end{gather*}
where $\{\SA^l\hat{\dot\x},\hat{\dot\P},\SA^l\check{\dot\x},\check{\dot\P}\}$ are given by a Kalman filter. Such a Kalman filter is then run for each hypothesis (sample) of $\SA^n\dot\x$ and $\SA^n\ddot\x$ in a particle filter.

The marginalized particle filters also handle special cases between assumption \emph{a} and \emph{c}. In case $\SA^n=\TA$, i.e. the output is independent of the linear subspace of the state equation, then the general output model of \emph{a} can be combined with the state equation of \emph{c}. \emph{Vice versa}, if $\SA=\TA^n$, i.e. the state equation is independent of the linear subspace of the output model, then the general state equation of \emph{a} can be combined with the output model of \emph{c}.

\subsection{Predecessors and applications}

A two step variant of the inactive subspace update marginalization for a special case of $\TA$ and assumptions $\mathcal{N}$ and \emph{a}, in terms of mean and quadratic mean, can be found in~\cite{Zachariah2012}. A variant for an arbitrary $\TA$ can be found in~\cite{Nilsson2013}. A one-dimensional variant, in terms of mean and covariance, can be found in~\cite{Nilsson2014}. In~\cite{Rantakokko2016}, the presented techniques are used for robust tightly-coupled GPS-aided dead reckoning.


\section{Conclusion}\label{sec:Conclusions}
In this article I have presented Bayesian filtering techniques which, given Gaussian priors and posteriors, enable one to exploit a set of sequentially more constraining model assumptions, condensing the filtering problem to essential moment integrals and enabling introduction of model approximations in the same integrals. Combined with suitable moment integral approximation techniques, e.g. quadrature or importance sampling, these techniques provide a new class of filtering techniques. In particular, in contrast to Kalman filters, the techniques can handle non-Gaussian (e.g. heavy-tailed) likelihood functions and, in contrast to marginalized particle filters, the techniques can handle dynamic subspace structures (e.g. state equations and output being non-linear in different subspaces). Further, the model assumptions enable system composition (combination of different submodels) by specifying subspace projection matrices, decoupling the \mbox{(sub-)} system modelling and the system composition. Moreover, the composition becomes transparent, enabling the Bayesian filtering techniques to exploit the underlying structure.


In summary, I have shown how model structure and model approximations, combined with marginalization, can reduce the filtering problem to low-dimensional moment integrals which enables us to handle a more general class of filtering problems as compared to Kalman filters and which open up for a large range of numerical approximation techniques. Further, for specific models and distributions, analytical solutions to the marginal moment integrals may be possible.



{\footnotesize
\bibliographystyle{ieeetr}
\bibliography{MSKF}
}

\end{document}